\documentclass[12pt,a4paper,fleqn]{article}

\usepackage{array}
\usepackage{natbib}
\usepackage{mathtools}
\renewcommand{\cite}{\citet}
\usepackage{float}

\usepackage{amsthm,amsmath}
\usepackage{bbm}
\usepackage{amssymb,subfig}
\usepackage{relsize}
\usepackage{graphicx}
\usepackage{epstopdf}
\usepackage{bigints}
\usepackage{url}
\allowdisplaybreaks

\usepackage{fancyhdr,enumerate}
\usepackage{verbatim}
\usepackage{latexsym}
\usepackage{amsfonts}
\usepackage[nohead]{geometry}
\usepackage[neverdecrease]{paralist}

\usepackage{multirow}

\usepackage{setspace}
\usepackage[bottom]{footmisc}
\usepackage{indentfirst}
\usepackage{endnotes}
\usepackage{rotating}
\usepackage{verbatim}
\usepackage{setspace}
\usepackage{multirow}
\usepackage{latexsym}
\usepackage{endnotes}
\RequirePackage[utf8x]{inputenc}

\usepackage{color}
\definecolor{dgreen}{rgb}{0,0.5,0}
\definecolor{dblue}{rgb}{0,0,0.9}
\definecolor{dred}{rgb}{0.6,0.0,0.1}
\definecolor{dgold}{rgb}{0.5,0.3,0.0}
\definecolor{dvio}{rgb}{0.6,0.3,0.5}
\definecolor{gray}{rgb}{0.5,0.5,0.5}

\newtheorem{hypo}{Assumption}[section]
\newtheorem{lem}{Lemma}[section]

\newtheorem{tr}{Theorem}[section]

\theoremstyle{definition}
\newtheorem{example}{Example}[section]
\newtheorem{rem}{Remark}[section]
\newtheorem{restr}{Restriction}

\newcommand{\la}{\langle}
\newcommand{\ra}{\rangle}
\newcommand{\E}{\mathbf{E}}
\newcommand{\1}{\mathbbm{1}}
\newtheorem{assump}{Assumption}

\usepackage{xr}
\numberwithin{equation}{section}

\begin{document}
\vspace{-3cm}
\title{Gaussian processes and Bayesian moment estimation\footnote{First version: February 2012. The authors gratefully thank the Editor, an Associate Editor, and three anonymous referees for their many constructive
comments on the previous version of the paper. The authors are grateful to Yuichi Kitamura, Frank Kleibergen, Andriy Norets and seminars and conferences participants at: Berlin, Boston College, Bristol, Carlos III, CREST, Leiden, Northwestern, SBIES 2015, ICEEE 2015, CFE - CMStatistics 2015, NASM 2012, Toulouse, for useful comments. We thank financial support from ANR-13-BSH1-0004 (IPANEMA). Anna Simoni gratefully acknowledges financial support from Labex ECODEC (ANR - 11-LABEX-0047), SFB-884 and hospitality from the University of Mannheim.}
}

\author{\begin{tabular}{ccc}
Jean-Pierre Florens\footnote{Toulouse School of Economics, Universit\'{e} de Toulouse Capitole, Toulouse - 21, all\'{e}e de Brienne - 31000 Toulouse (France). Email: jean-pierre.florens@tse-fr.eu} & & Anna Simoni\footnote{CREST, CNRS, \'{E}cole Polytechnique, ENSAE - 5, avenue Henry Le Chatelier, 91120 Palaiseau, France. Email: simoni.anna@gmail.com (\textit{corresponding author}).}
\end{tabular}}

\date{}
\maketitle

\begin{abstract}
Given a set of moment restrictions (MRs) that overidentify a parameter $\theta$, we investigate a semiparametric Bayesian approach for inference on $\theta$ that does not restrict the data distribution $F$ apart from the MRs. As main contribution, we construct a degenerate Gaussian process prior that, conditionally on $\theta$, restricts the $F$ generated by this prior to satisfy the MRs with probability one. Our prior works even in the more involved case where the number of MRs is larger than the dimension of $\theta$. We demonstrate that the corresponding posterior for $\theta$ is computationally convenient. Moreover, we show that there exists a link between our procedure, the Generalized Empirical Likelihood with quadratic criterion and the limited information likelihood-based procedures. We provide a frequentist validation of our procedure by showing consistency and asymptotic normality of the posterior distribution of $\theta$. The finite sample properties of our method are illustrated through Monte Carlo experiments and we provide an application to demand estimation in the airline market.\\
%
\end{abstract}

\section{Introduction}
\indent Econometric models are often formulated via moment restrictions (MRs) of the form $\mathbf{E}^F[h(\theta,X)] = 0$, where $h(\theta,X)$ is a vector-valued function of an observable random element $X$ with distribution $F$ and a parameter vector $\theta$. These MRs provide the only information available about $\theta$ and the data distribution. Given a set of MRs, this paper builds a Bayesian inference procedure for $\theta$ that imposes these moment conditions in the nonparametric prior distribution for the data distribution $F$ and that is computationally convenient. Apart from these MRs, $F$ is left unrestricted.\\
\indent A main advantage of Bayesian inference consists of providing a well-defined posterior distribution which is important for many decision problems. On the other hand, constructing Bayesian inference procedures for moment condition models presents two difficulties. A first difficulty is due to the fact that a likelihood is not available. A second difficulty arises because imposing overidentifying MRs on the prior distribution for the nonparametric $F$ is challenging. The contribution of this paper is to propose an elegant approach that allows to deal with these two difficulties.\\ 
\indent The model we consider is as follows. Let $X$ be an observable random element in $\mathbb{R}^{m}$ with distribution $F$ and $X_{1},\ldots,X_{n}$ be an \textit{i.i.d.} sample of $X$. The parameter $\theta\in\Theta\subset \mathbb{R}^p$ is linked to the data generating process (DGP) $F$ through the MRs
    \begin{equation}\label{eq_moment_condition_introduction}
      \mathbf{E}^{F}\left[h(\theta,X)\right] = 0,
    \end{equation}
\noindent where $h(\theta,x) := (h_1(\theta,x),\ldots,h_d(\theta,x))^T$ and the functions $h_j(\theta,x)$, $j=1,\ldots,d$ are real-valued and known up to $\theta$. We assume $d\geq p$ and our main interest is the case where $d>p$, which is in general more challenging than the case $d=p$, see \textit{e.g.} in \cite{BroniatowskiKeziou2012}. Apart from \eqref{eq_moment_condition_introduction}, $F$ is completely unrestricted.\\ 
\indent Imposing MRs via semiparametric priors may be challenging depending on the relationship existing between $\theta$ and $F$. More precisely, when the model is exactly identified (\textit{i.e.} $p = d$) and \eqref{eq_moment_condition_introduction} characterizes $\theta$ as an explicit function of $F$, say $\theta = b(F)$ for some function $b(\cdot)$, then \eqref{eq_moment_condition_introduction} does not restrict $F$. Therefore, one places an unrestricted nonparametric prior on $F$ and recovers the prior of $\theta$ via the transformation $\theta = b(F)$. The $(\theta,F)$s generated by this prior automatically satisfy the constraints.\\
\indent On the contrary, when $d>p$ (overidentified model), $\theta$ cannot be expressed as an explicit function of $F$. Indeed, \eqref{eq_moment_condition_introduction} imposes constraints on $F$ and the existence of a solution $\theta$ to \eqref{eq_moment_condition_introduction} is guaranteed only for a subset of distributions $F$. Therefore, a restricted nonparametric prior on $F$ must be specified conditionally on $\theta$ and the support of this prior is a proper subset of the set of probability distributions. It turns out that incorporating overidentifying MRs in a semiparametric prior for $(\theta,F)$ is not straightforward. In this paper we propose a way to construct a semiparametric prior that incorporates the overidentifying MRs. With this prior we show that the corresponding marginal posterior of $\theta$ is proportional to the following expression:
\begin{multline}
  \exp\Big\{ - \frac{n}{2}\Big[\frac{1}{n}\sum_{i=1}^n[h(\theta,x_i)^T] \left(\E^*[h(\theta,X)h(\theta,X)^T]\right)^{-1}\frac{1}{n}\sum_{i=1}^n[h(\theta,x_i)]\\
  + \sum_{j>d}\left(\frac{1}{n}\sum_{i=1}^{n}\varphi_j(x_i) - \mathbf{E}^*[\varphi(f_{0\theta})(X)\varphi_j(X)]\right)^2\frac{1}{1 + n\lambda_j} - \frac{2\log\mu(\theta)}{n}\Big]\Big\}\label{eq_posterior_introduction}
\end{multline}
where $\mathbf{E}^*$ denotes the expectation taken with respect to the true data distribution, $(\varphi_j)_{j>d}$ are functions orthogonal to the span of $(1,h(\theta,X)^T)^T$, and $\varphi(f_{0\theta})$, $\lambda_j$ and $\mu(\theta)$ are quantities coming from our semiparametric prior that will be defined below. There are many remarkable facts associated with this expression. (I) The first term in the exponential is proportional to the continuous-updating GMM criterion, and for a fixed prior, the last two terms vanish as $n\rightarrow \infty$. (II) The second term in the exponential accounts for the extra information contained in the subspace orthogonal to the subspace spanned by $(1,h(\theta,X)^T)^T$. This extra information is brought by the prior and can be exploited even asymptotically if the prior is not fixed but varies with $n$ at an appropriate rate. (III) Expression \eqref{eq_posterior_introduction} makes a connection between the parametric case (when the second term does not vanish asymptotically) and the semiparametric case (when the second term vanishes). (IV) Expression \eqref{eq_posterior_introduction} is particularly convenient when one wants to sample from the posterior distribution by using a Metropolis-Hasting algorithm.\\
\indent While the mathematical construction that we develop to obtain \eqref{eq_posterior_introduction} is technically involved, one only needs to use \eqref{eq_posterior_introduction} to implement our approach. Our mathematical construction is based on a degenerate Gaussian process ($\mathcal{GP}$) prior with restricted support which is easy to deal with and that works as follows. The DGP $F$ is assumed to admit a density function $f$ with respect to some positive measure $\Pi$ chosen by the researcher (for instance the Lebesgue measure). Then, we endow $f$ with a $\mathcal{GP}$ prior conditional on $\theta$. The $d\geq p$ MRs are incorporated by constraining the prior mean and prior covariance of this $\mathcal{GP}$ in an appropriate way. Because this prior imposes the MRs, it is degenerate on a proper subset of the set of functions integrating to one. The reason for the appropriateness of a $\mathcal{GP}$ prior in such a framework is due to the fact that the MRs in \eqref{eq_moment_condition_introduction} are linear in $f$ and the linearity of the model matches extremely well with a $\mathcal{GP}$ prior. A remarkable feature of our method is that the MRs are imposed directly through the $\mathcal{GP}$ prior of $f$ given $\theta$ without requiring a second step projection over the set of density functions satisfying the MRs. To the best of our knowledge a $\mathcal{GP}$ prior has not been used yet in the MRs framework.\\
\indent Our Bayesian procedure, that we call the $\mathcal{GP}$-approach, is constructed as follows. We first specify a prior on $\theta$ and then a $\mathcal{GP}$ prior on $f$ conditional on $\theta$. We circumvent the difficulty of the likelihood function specification, which is not available, by constructing a linear functional transformation of the DGP $F$ such that its empirical counterpart, say $r_n$, has an asymptotic Gaussian distribution. This will be used as the sampling model. Therefore, our model is approximately conjugate and allows easy computations while being nonparametric in $F$.\\
\indent We provide a closed-form expression for the marginal posterior distribution of $\theta$, obtained by integrating out $f$ and proportional to \eqref{eq_posterior_introduction}. We show that the quasi-likelihood function for GMM settings -- which is the finite sample analog of the limited information likelihood (LIL) -- used, among others, by \citet{Kim2002} and \citet{ChernozhukovHong2003}, can be obtained as the limit of our marginal likelihood function when the $\mathcal{GP}$ prior for $f$ is allowed to become diffuse. In addition, when the prior for $f$ becomes noninformative, the marginal posterior distribution for $\theta$ becomes the same (up to constants) as the GEL objective function with quadratic criterion and is a monotonic transformation of the continuous-updating GMM objective function (\citet{HansenHeatonYaron1996}).\\
\indent Finally, we provide a frequentist validation of our method by showing: \textit{(i)} posterior consistency, \textit{(ii)}  frequentist consistency of the maximum a posteriori estimator, and \textit{(iii)} asymptotic normality of the posterior distribution of $\theta$.
\paragraph{Related literature.}
Inference in an MRs framework has received very much attention in the past literature. Among the most popular statistical methods in this framework is the Empirical Likelihood (EL)-based inference (\textit{e.g.} \citet{Owen1988}, \citet{QinLawless1994}, \cite{Imbens1997}, \citet{Owen2001}, \citet{ChenVanKeilegom2009}) and its variations  like the Exponential Tilting (\textit{e.g.} \citet{KitamuraStutzer1997} and \citet{Kitamura1997}), the Exponentially Tilted EL (ETEL) proposed by \citet{Schennach2007}, and the Generalized EL (GEL) (\textit{e.g.} \citet{Smith1997}, \citet{NeweySmith2004}). EL and ETEL can also be used in a Bayesian framework for valid inference as established by \citet{Lazar2003}, \citet{grendar2009}, \citet{Schennach2005}, \citet{Rao2010}, \citet{ChibShinSimoni2017} and \citet{Chaudhuri2017} among others.\\
\indent An alternative frequentist inference procedure, which is very popular especially in econometrics, is the Generalized method of Moments (GMM) (\citet{Hansen1982} and \citet{HansenSingleton1982}). From a Bayesian inference point of view, the semiparametric procedures that have been proposed in the literature can be classified in two types depending on whether the moment conditions \eqref{eq_moment_condition_introduction} are imposed in the likelihood or in the nonparametric prior. Procedures of the first type construct a quasi-likelihood by exponentiating either the quadratic criterion or the GEL criterion associated with the empirical counterpart of \eqref{eq_moment_condition_introduction} and include \citet{Kwan1999}, \citet{Kim2002}, \citet{ChernozhukovHong2003} and \citet{liao2011} among others. Our paper shows that the quasi-likelihood used in this type of approach arises as the limit of our $\mathcal{GP}$ prior as it becomes diffuse. We provide thus a fully Bayesian justification to this approach. Procedures of the second type impose the moment conditions in the prior for $(\theta,F)$ while leaving the likelihood completely unrestricted and include, in addition to the Bayesian EL and ETEL discussed above, \citet{ChamberlainImbens2003} who use a Dirichlet prior, \citet{KitamuraOtsu2011} and \citet{Shin2014} who propose a two-step procedure based on a projection of a Dirichlet process mixture prior and of a mixture of Dirichlet Process prior, respectively, and \citet{BornnShephardSolgi2015} who use Hausdorff measures to build probability tools for dealing with moment estimation. Parametric Bayesian procedure have been proposed for the linear instrumental variable regression model and were based on a parametric specification of the error distributions, see \textit{e.g.} \cite{KleibergenZivot2003} and references therein.\\
\indent The paper is organized as follows. The $\mathcal{GP}$-approach is described in sections \ref{s_general_model} and \ref{s_Posterior_Distribution}, which contain our main theoretical contribution. In section \ref{s_Posterior_Distribution} we also show the link existing between our approach and some frequentist approaches in an MRs framework, as mentioned above. In section \ref{s_asymptotic_analysis} we analyze frequentist asymptotic properties of the posterior distribution of $\theta$. In section \ref{s_implementation} we provide an MCMC algorithm to implement our method, develop simulation studies and provide an application to demand estimation in the airline market. Section \ref{s_conclusions} concludes. A testing procedures, further details of our $\mathcal{GP}$-approach and all the proofs are gathered in the Supplementary Material.
%
\section{The Gaussian Process ($\mathcal{GP}$)-approach}\label{s_general_model}
\indent Let $X$ be a continuous random element in $S\subseteq\mathbb{R}^{m}$ with distribution $F$ and $X_{1}, \ldots,X_{n}$ be an i.i.d. sample of $X$. Assume that $F$ is absolutely continuous with respect to some positive measure $\Pi$ (\textit{e.g.} the Lebesgue measure) with density function $f$. In other words, conditionally on $f$ the data are drawn from $F$: $X_{1}, \ldots,X_{n}|f\sim F$. The set of probability density functions (\textit{pdf}s) on $S$ with respect to $\Pi$ is denoted by $M$.\\
\indent Let $\theta\in\Theta\subseteq\mathbb{R}^p$ be the parameter of interest characterized by (\ref{eq_moment_condition_introduction}). By adopting a frequentist point of view, we denote, throughout the paper, the true value of $\theta$ by $\theta_{*}$, the true DGP by $F_{*}$ and its density with respect to $\Pi$ by $f_{*}$. The model is assumed to be well-specified, that is, there exists $\theta_*\in\Theta$ such that $\mathbf{E}^{F_{*}}(h(\theta_{*},X)) = 0$ holds. We endow $S\subseteq\mathbb{R}^{m}$ with the trace of the Borelian $\sigma$-field $\mathfrak{B}_{S}$ and specify $\Pi$ as a positive measure on this subset. We denote by $\mathcal{E} := L^{2}(S,\mathfrak{B}_{S},\Pi)$ the Hilbert space of square integrable functions on $S$ with respect to $\Pi$ and by $\mathfrak{B}_{\mathcal{E}}$ the Borel $\sigma$-field generated by the open sets of $\mathcal{E}$. The scalar product and norm on this space are defined in the usual way and denoted by $\langle\cdot,\cdot\rangle$ and $||\cdot||$, respectively.\\
\indent The parameters of the model are $(\theta,f)$, where $f$ is the nuisance parameter, and the parameter space is $\Lambda := \left\{(\theta, f)\in\Theta\times \mathcal{E}_{M}; \int h(\theta, x)f(x)\Pi(dx) = 0\right\}$, $\mathcal{E}_{M} := \mathcal{E}\cap M$, where $h:\Theta\times \mathbb{R}^{m}\rightarrow \mathbb{R}^{d}$ is a known function. In the following of the paper we maintain the following assumption.
    \begin{hypo}\label{Ass_2_1}
      \textit{(i)} The true $f_{*}$ satisfies $f_{*}\in\mathcal{E}_{M}:=\mathcal{E}\cap M$; \textit{(ii)} the moment function $h(\theta,\cdot)$ is such that $h_{i}(\theta,\cdot)\in \mathcal{E}$ for every $i=1,\ldots,d$ and for every $\theta\in\Theta$, where $h_{i}$ denotes the $i$-th component of $h$; \textit{(iii)} $d\geq p$.
    \end{hypo}
\indent Assumption \ref{Ass_2_1} \textit{(i)} restricts $f_{*}$ to be square integrable with respect to $\Pi$ and is for instance verified if $f_{*}$ is bounded and $\Pi$ is a bounded measure. The model is made up of three elements that we detail in the next two subsections: a prior on $\theta$, denoted by $\mu(\theta)$, a conditional prior on $f$ given $\theta$, denoted by $\mu(f|\theta)$, and the sampling model. In the following, we shorten ``almost surely'' by ``a.s.'' and omit the probability which ``a.s.'' refers to. We denote by $\mathbf{E}^F$ the expectation taken with respect to $F$ and by $\mathbf{E}^*$ the expectation taken with respect to $F_*$.
%
\subsection{Prior distribution}\label{ss_prior_distribution}
We specify a prior probability measure $\mu$ for $(\theta,f)$ of the form $\mu(\theta,f) = \mu(\theta) \mu(f|\theta)$. By abuse of notation, $\mu(\theta)$ will also denote the Lebesgue density of the prior distribution of $\theta$ in the case it admits it. The prior $\mu(\theta)$ may either be flat (non-informative) or incorporate any additional information available to the econometrician about $\theta$. In any case, it is tacitly assumed that $\mu(\theta)$ is such that the posterior of $\theta$ exists.\\
\indent Given a value for $\theta$, the conditional prior $\mu(f|\theta)$ is specified such that its support equals the subset of functions in $\mathcal{E}$ that integrate to one and satisfy (\ref{eq_moment_condition_introduction}) for this particular value of $\theta$. At the best of our knowledge, the construction of such a conditional prior $\mu(f|\theta)$ is new in the literature and we now explain it in detail.
\paragraph{Construction of the conditional prior $\mu(f|\theta)$.} We construct the conditional prior distribution $\mu(f|\theta)$ of $f$, given $\theta$, as a $\mathcal{GP}$ on $\mathfrak{B}_{\mathcal{E}}$ with mean function $f_{0\theta}\in\mathcal{E}_{M}$ and covariance operator $\Omega_{0\theta}:\mathcal{E}\rightarrow \mathcal{E}$. We restrict $f_{0\theta}$ and $\Omega_{0\theta}$ to guarantee that the trajectories $f$ generated by $\mu(f|\theta)$ are such that the corresponding $F$ (which is such that $dF = f d\Pi$) integrates to $1$ and satisfies equation \eqref{eq_moment_condition_introduction} with probability $1$. The two sets of restrictions that we impose are the following (one on $f_{0\theta}$ and one on $\Omega_{0\theta}$):
\begin{restr}[Restriction on $f_{0\theta}$]\label{R1}
  The prior mean function $f_{0\theta}\in\mathcal{E}_{M}$ is chosen such that
    \begin{equation}\label{eq_moment_condition}
      \int h(\theta,x)f_{0\theta}(x)\Pi(dx) = 0.
    \end{equation}
\end{restr}
\begin{restr}[Restriction on $\Omega_{0\theta}$]\label{R2}
  The prior covariance operator $\Omega_{0\theta}:\mathcal{E}\rightarrow\mathcal{E}$ is chosen such that
    \begin{equation}\label{eq_prior_covariance_restrictions}
      \left\{\begin{array}{ccc}
        \Omega_{0\theta}^{1/2}h(\theta,x) & = & 0\\
        \Omega_{0\theta}^{1/2}1 & = & 0
      \end{array}\right.
    \end{equation}
  where $\Omega_{0\theta}^{1/2}:\mathcal{E}\rightarrow\mathcal{E}$ denotes the positive square root of $\Omega_{0\theta}$: $\Omega_{0\theta} = \Omega_{0\theta}^{1/2}\Omega_{0\theta}^{1/2}$.
\end{restr}
The covariance operator $\Omega_{0\theta}$ is linear, self-adjoint and trace-class.\footnote{A trace-class operator is a compact operator with eigenvalues that are summable. Remark that this guarantees that the trajectories $f$ generated by $\mu(f|\theta)$ satisfy $\int f^{2}d\Pi<\infty$ a.s.} Due to Restriction \ref{R2}, $\Omega_{0\theta}$ is not injective. In fact, the null space of $\Omega_{0\theta}$, denoted by $\mathfrak{N}(\Omega_{0\theta})$, is not trivial and contains effectively the constant $1$ -- which implies that the trajectories $f$ generated by the prior integrate to $1$ a.s. (with respect to $\Pi$) -- and the functions $h(\theta,x)$ -- which implies that the trajectories $f$ satisfy the moment conditions a.s. This means that $\Omega_{0\theta}$ is degenerate in the directions along which we want that the corresponding projections of $f$ and $f_{0\theta}$ are equal. Therefore, the support of $\mu(f|\theta)$ is a proper subset of $\mathcal{E}$. This is stated in the next lemma.
\begin{lem}\label{lem_prior}
  The conditional $\mathcal{GP}$ prior $\mu(f|\theta)$, with mean function $f_{0\theta}$ and covariance operator $\Omega_{0\theta}$ satisfying Restrictions \ref{R1} and \ref{R2}, generates trajectories $f$ that satisfy $\mu(f|\theta)$-a.s. the conditions $\int f(x)\Pi(dx) = 1$ and $\int h(\theta,x)f(x)\Pi(dx) = 0$, for every $\theta\in\Theta$.
\end{lem}
\begin{rem}
  Restrictions \ref{R1} and \ref{R2} imply that the trajectories generated by $\mu(f|\theta)$ integrates to $1$ (with respect to $\Pi$) and satisfy \eqref{eq_moment_condition_introduction} a.s. but they do not guarantee non-negativity of the trajectories. To impose non-negativity one could: (i) either project the restricted $\mathcal{GP}$ prior on the space of non-negative functions, or (ii) write either $f = g^{2}$ or $f = e^g/\int e^g$, $g\in\mathcal{E}$, and specify a conditional $\mathcal{GP}$ prior distribution for $g$, given $\theta$, instead of for $f$. Nonetheless, it is important to notice that: the projected prior in (i) is not Gaussian anymore, and in (ii) we cannot use our restricted $\mathcal{GP}$ prior for $g$ because it does not work anymore to impose the MRs on $f$ so that a different restricted prior should be constructed. Moreover, in both (i) and (ii) the resulting posterior for $\theta$ is not available in closed form which is instead one of the main advantages of our procedure. Therefore, it is not possible to impose the non-negativity constraint if one wants to use our restricted $\mathcal{GP}$ prior. However, because our goal is to make inference on $\theta$ while $f$ is a nuisance parameter, failing to impose the non-negativity constraint is not an issue as long as our procedure is shown to be consistent for $\theta$ (which we show in section \ref{s_asymptotic_analysis}). Of course, the finite sample performance will be affected by our choice of the prior on $f$ (as well as on $\theta$). However, non reliability of finite sample performance is in general the case for inference procedures based on MR models like the GMM and the GEL.
\end{rem}
\indent From a practical implementation point of view, a covariance operator satisfying Restriction \ref{R2} and a $f_{0\theta}$ satisfying Restriction \ref{R1} may be constructed as follows.\\
\textbf{Construction of $\Omega_{0\theta}$.} Let $(\lambda_{j})_{j\in\mathbb{N}}$ be a decreasing sequence of non-negative numbers such that $\sum_j\lambda_j<\infty$, and $(\varphi_{j})_{j\in\mathbb{N}}$ be an orthonormal basis (o.n.b.) for $\mathcal{E}$. Then, $\forall \phi \in \mathcal{E}$: $\Omega_{0\theta} \phi = \sum_{j=0}^{\infty}\lambda_{j} \langle \phi,\varphi_{j}\rangle\varphi_{j}$. Remark that $(\lambda_{j})_{j\in\mathbb{N}}$ and $(\varphi_{j})_{j\in\mathbb{N}}$ correspond to the eigenvalues and eigenfunctions of $\Omega_{0\theta}$, respectively. Since the null space $\mathfrak{N}(\Omega_{0\theta})\subset\mathcal{E}$ is spanned by $\{1,h_{1}(\theta,\cdot),\ldots,h_{d}(\theta,\cdot)\}$, we can set the first eigenfunctions of $\Omega_{0\theta}$ equal to the elements of any o.n.b. of $\mathfrak{N}(\Omega_{0\theta})$. Restriction \ref{R2} is then fulfilled by setting the corresponding eigenvalues equal to $0$. For instance, if $\{1,h_{1}(\theta,\cdot),\ldots,h_{d}(\theta,\cdot)\}$ are orthonormal as elements of $\mathcal{E}$, then $\mathfrak{N}(\Omega_{0\theta})$ has dimension $d+1$, the first eigenfunctions are $(\varphi_{0},\varphi_{1}, \ldots, \varphi_{d})^{T} = (1,h^T)^T$ and the corresponding eigenvalues are $\lambda_{j} = 0$, $\forall j = 0,1,\ldots, d$. Remark that in this case, necessarily, $\int\Pi(dx) = 1$. If $\{1,h_{1}(\theta,\cdot),\ldots,h_{d}(\theta,\cdot)\}$ are not orthonormal then one can use their orthonormalized counterparts as the first eigenfunctions of $\Omega_{0\theta}$. The remaining components $(\varphi_{j})_{j>d}$ are chosen such that $(\varphi_{j})_{j\geq 0}$ forms an orthonormal basis of $\mathcal{E}$ and $(\lambda_{j})_{j>d}$ are chosen such that $\sum_{j>d}\lambda_{j} <\infty$. Hence, $\forall \phi \in\mathcal{E}$, we construct $\Omega_{0\theta}$ as $\Omega_{0\theta} \phi = \sum_{j=d+1}^{\infty}\lambda_{j} \langle \phi,\varphi_{j}\rangle\varphi_{j}$ where we suppress the dependence of $\varphi_j$ on $\theta$ for simplicity. Examples of choices for $(\lambda_{j})_{j>d}$ are, for some constant $c>0$: \textit{(i)} $\lambda_j = c j^{-a}$ with $a>1$, \textit{(ii)} $\lambda_j = ce^{-j}$. In section \ref{s_implementation} we give more details about the construction of $\Omega_{0\theta}$.\\
\textbf{Construction of $f_{0\theta}$.} When the MRs are multivariate polynomial functions of the components of $X$ then a simple way to construct $f_{0\theta}$ is to set it equal to a $\mathcal{N}(0,V)$ with some of the elements of $V$ restricted to satisfy the MRs. 
\begin{rem}\label{Remark_2.3}
  In addition to the parameters characterized by the MRs, the prior distribution of $f$ can depend on other hyperparameters, say $\beta$. In this case $\theta:=(\beta,\widetilde\theta)^T$ where $\widetilde\theta$ denotes the parameters characterized by the MRs. Our procedure does not change in this case and all the computations and results of the paper continue to hold.
\end{rem}
%
\subsection{A tractable sampling model}\label{ss_sampling_model}
Given the observed \textit{i.i.d.} sample $(x_1,\ldots,x_n)$, the likelihood function is $\prod_{i=1}^{n}f(x_{i})$. While apparently simple, using this likelihood for Bayesian inference on $\theta$ makes the analysis of the posterior distribution complicated. This is because to compute the posterior for $\theta$ one has to marginalize out $f$. Since a $\mathcal{GP}$ prior is not a natural conjugate of the \textit{i.i.d.} model then, marginalization of $f$ has to be carried out through numerical, or Monte Carlo, integration on a functional space, which may be computationally costly. To avoid this difficulty, we propose an alternative and original way to construct the sampling model that allows for a conjugate analysis and prevents from numerical integration. Our approach is based on a functional transformation $r_n$ of the sample $X_{1},\ldots, X_{n}$. This sampling model allows us to obtain the tractable posterior of $\theta$ proportional to \eqref{eq_posterior_introduction}.\\
\indent The transformation $r_n$ is chosen by the researcher and must have the following characteristics: \textit{(I)} $r_n$ is an observable element of an infinite-dimensional Hilbert space $\mathcal{F}$ (to be defined below), for instance a $L^{2}$-space; \textit{(II)} $r_n$ converges weakly towards a Gaussian process in $\mathcal{F}$; \textit{(III)} the expectation of $r_n$, conditional on $f$, defines a linear operator $K:\mathcal{E}\rightarrow\mathcal{F}$ such that $\mathbf{E}^F(r_n) = Kf$; \textit{(IV)} $r_n$ is a one-to-one transformation of some sufficient statistic. Moreover, $r_n\in\mathcal{F}$ is a Hilbert space-valued random variable (H-r.v.), that is, for a complete probability space $(Z,\mathcal{Z},\mathbb{P})$, $r_n$ defines a measurable map $r_n:(Z,\mathcal{Z},\mathbb{P})\rightarrow(\mathcal{F},\mathfrak{B}_{\mathcal{F}})$, where $\mathfrak{B}_{\mathcal{F}}$ denotes the Borel $\sigma$-field generated by the open sets of $\mathcal{F}$.
\paragraph{Construction of $r_n$.} Let $\mathfrak{T}\subseteq \mathbb{R}^{l}$, $l>0$. To construct $r_n$ we first select a function $k(t,x): \mathfrak{T}\times S\rightarrow \mathbb{R}$ (or in $\mathbb{C}$) that is measurable in $x$ for every $t\in \mathfrak{T}$ and that is non-constant in $(t,x)$. The transformation $r_n$ is then taken to be the expectation of $k(t,\cdot)$ under the empirical measure:
    \begin{displaymath}
      r_n(t) := \frac{1}{n}\sum_{i=1}^{n}k(t,x_{i}), \qquad \forall t\in\mathfrak{T}.
    \end{displaymath}
\noindent Define $\mathcal{F} = L^{2}(\mathfrak{T},\mathfrak{B}_{\mathfrak{T}},\rho)$ where $\rho$ is a measure on $\mathfrak{T}$ and $\mathfrak{B}_{\mathfrak{T}}$ denotes the Borel $\sigma$-field generated by the open sets of $\mathfrak{T}$. The scalar product and norm on $\mathcal{F}$ are defined in the usual way and denoted by $\langle\cdot,\cdot\rangle$ and $\|\cdot\|$, respectively, by using the same notation as for the inner product and norm in $\mathcal{E}$. The function $k(t,\cdot)$ defines also a bounded operator $K$: 
\begin{equation}\label{eq_K_definition}
  \begin{array}{rcl}
    K:\mathcal{E} & \rightarrow & \mathcal{F}\\
    \varphi & \mapsto & \int k(t,x)\varphi(x)\Pi(dx)
  \end{array}
\end{equation}
\noindent and must be such that, for every $\varphi\in\mathcal{E}$, $K\varphi\in\mathcal{F}$ and $r_n$ is an H-r.v. with realizations in $\mathcal{F}$. For every $f\in\mathcal{E}_M$, $Kf$ is the expectation of $k(t,X)$ under $F$: $(Kf)(t) = \mathbf{E}^{F}(k(t,X))$. Under the true distribution $F^*$ the expectation of $r_n$ is $Kf_*$ and the covariance function of $r_n$ is: $\forall s,t\in\mathfrak{T}$, $\frac{1}{n}\sigma(t,s) := \mathbf{E}^{*}[r_n(t) - (Kf_*)(t)][r_n(s) - (Kf)(s)]$  $= \frac{1}{n}\mathbf{E}^{*}\left[k(t,X)k(s,X)\right] - \frac{1}{n}\mathbf{E}^{*}[k(t,X)]\mathbf{E}^{*}[k(s,X)].$
\noindent If the class of functions $\{k(t,\cdot),t\in\mathfrak{T}\}$ is Donsker then, as $n\rightarrow \infty$, the conditional distribution of $\sqrt{n}(r_n - Kf_*)$ weakly converges to a $\mathcal{GP}$ with covariance operator $\Sigma:\mathcal{F}\rightarrow \mathcal{F}$ defined as
    \begin{equation}
      \forall \psi\in\mathcal{F},\qquad (\Sigma \psi)(t) = \int \sigma(t,s)\psi(s)\rho(ds)\label{approximated_sampling_model}
    \end{equation}
\noindent which is one-to-one, linear, positive definite, self-adjoint and trace-class. In the following, we assume that $\{k(t,\cdot),t\in\mathfrak{T}\}$ is Donsker so that $r_n$ is approximately Gaussian: $r_n \sim \mathcal{GP}(Kf_*,\Sigma_n)$ where $\Sigma_n = \frac{1}{n}\Sigma$. Finally, among all the functions $k(t,x)$ that satisfy the previous assumptions, one should keep only the ones such that the corresponding $r_n$ is a one-to-one transformation of some sufficient statistic, as required in \textit{(IV)} above. This, which will be tacitly assumed in the following of the paper, guarantees that the posterior distribution computed by using $r_n$ does not depend on the particular choice of $k(t,x)$.\\
\indent In our analysis we treat $f_*$ as a realization of the random parameter $f$ and $\Sigma_n$ as known. Therefore, the sampling distribution of $r_n|f$ is $P^f = \mathcal{GP}(Kf,\Sigma_n)$ and we construct the posterior distribution based on it. In practice, $\Sigma_n$ must be replaced by its empirical counterpart. In finite sample, $P^f$ is an approximation of the true sampling distribution but the approximation error vanishes as $n\rightarrow \infty$. Moreover, the approximating sampling distribution $P^f$ is only used to construct the posterior distribution and the proofs of our asymptotic results do not rely on it. We give now two examples where $r_n$ is a one-to-one transformation of sufficient statistics and $K$ is injective.

\begin{example}[Empirical cumulative distribution function (cdf)]
  Let $(X_{1}, \ldots, X_{n})$ be an \textit{i.i.d.} sample of $X\in\mathbb{R}$. A possible choice for $k(t,x)$ is $k(t,x) = 1\{x\leq t\}$, where $1\{A\}$ denotes the indicator function of the event $A$. In this case, $r_n(t) = F_n(t):=\frac{1}{n}\sum_{i=1}^{n}1\{X_{i}\leq t\}$ is the empirical \textit{cdf} and the operator $K$ is $(K\varphi)(t)= \int_{S} 1\{s\leq t\}\varphi(s)\Pi(ds)$, $\forall \varphi\in\mathcal{E}$. By the Donsker's theorem, $F_n(\cdot)$ is asymptotically Gaussian with mean the true cdf $F_*(\cdot)$ and covariance operator characterized by the kernel: $\frac{1}{n}(F_*(s\wedge t) - F_*(s)F_*(t))$.
\end{example}
\begin{example}[Empirical characteristic function]
  Let $(x_{1}, \ldots, x_{n})$ be an \textit{i.i.d.} sample of $x\in\mathbb{R}$. Let $k(t,x) = e^{itx}$, so that $r_n(t) = c_n(t) := \frac{1}{n}\sum_{j=1}^{n}e^{itx_{j}}$ is the empirical characteristic function. In this case, the operator $K$ is $(K\varphi)(t) = \int_{S} e^{its}\varphi(s)\Pi(ds)$, $\forall \varphi\in\mathcal{E}$. By the Donsker's theorem, $c_n(\cdot)$ is asymptotically a Gaussian process with mean the true characteristic function $c(\cdot) := \mathbf{E}^*[e^{itx}]$ and covariance operator characterized by the kernel: $\frac{1}{n}(c(s + t) - c(s)c(t))$.
\end{example}

\section{Posterior distribution}\label{s_Posterior_Distribution}
The Bayesian model defines a joint distribution on $(\theta,f,r_n)$ and can be summarized in the following way: $\theta\sim \mu(\theta)$, $f|\theta \sim \mu(f|\theta) = \mathcal{GP}(f_{0\theta},\Omega_{0\theta}),$ such that $\int h(\theta,x)f_{0\theta}(x)\Pi(dx) = 0$ and $\Omega_{0\theta}^{\frac{1}{2}}(1, h(\theta,\cdot)^T)^T = 0$, $r_n|f, \theta \sim r_n|f \sim P^{f} = \mathcal{GP}(Kf,\Sigma_{n})$ where we use the $\mathcal{GP}$ approximation $P^{f}$. Theorem 1 in \citet{FS2012SJS} shows that the joint distribution of $(f,r_n)$, conditional on $\theta$, is:
    \begin{equation}\label{eq_prior_distribution}
      \left. \left(\begin{array}
        {c} f\\ r_n
      \end{array}\right)\right|\theta \; \sim\; \mathcal{GP}\left(\left(\begin{array}
        {c}f_{0\theta}\\ K f_{0\theta}
      \end{array}\right), \left(\begin{array}
        {cc}\Omega_{0\theta}, & \Omega_{0\theta}K^{*}\\
        K\Omega_{0\theta}, & \Sigma_{n} + K\Omega_{0\theta}K^{*}
      \end{array}\right)\right)
    \end{equation}
\noindent where $(\Sigma_{n} + K\Omega_{0\theta}K^{*}):\mathcal{F}\rightarrow \mathcal{F}$, $\Omega_{0\theta}K^{*}: \mathcal{F}\rightarrow\mathcal{E}$, $K\Omega_{0\theta}:\mathcal{E}\rightarrow\mathcal{F}$ and $K^{*}:\mathcal{F} \rightarrow \mathcal{E}$ is the adjoint operator of $K$. We recall that the adjoint $K^{*}$ of a bounded and linear operator $K:\mathcal{E}\rightarrow \mathcal{F}$ is defined as the operator from $\mathcal{F}$ to $\mathcal{E}$ that satisfies $\langle K\varphi,\psi\rangle = \langle\varphi,K^{*}\psi\rangle$, $\forall \varphi\in\mathcal{E}$ and $\forall\psi\in\mathcal{F}$. In our case, $(K^{*}\psi)(t) = \int_{\mathfrak{T}} k(t,x)\psi(t)\rho(dt)$, $\forall\psi\in\mathcal{F}$. The marginal sampling distribution of $r_n$ conditional on $\theta$, obtained by integrating out $f$, is:
    \begin{equation}\label{eq_marginal_sampling_distribution}
      r_n|\theta\sim P_n^{\theta} := \mathcal{GP}(Kf_{0\theta}, \Sigma_{n} + K\Omega_{0\theta}K^{*}).
    \end{equation}
\noindent The conditional posterior distribution $\mu(f|r_n,\theta)$ of the nuisance parameter $f$ is briefly discussed in Appendix I and its properties follow from \cite{FS14}.\\
\indent The marginal posterior for $\theta$, denoted by $\mu(\theta|r_n)$ and which is the distribution of interest, is obtained by using the marginal sampling distribution $P_{n}^{\theta}$ given in (\ref{eq_marginal_sampling_distribution}). We first have to characterize the likelihood of $P_n^{\theta}$ with respect to an appropriate common dominating measure. This is done in Theorem A.1 in Appendix A and we denote by $p_{n\theta}(r_n;\theta)$ the likelihood function characterized in this theorem and used to construct $\mu(\theta|r_n)$. Hence, the marginal posterior distribution of $\theta$ writes:
$\mu(\theta|r_n) = \frac{p_{n\theta}(r_n;\theta) \mu(\theta)}{\int_{\Theta}p_{n\theta}(r_n;\theta) \mu(\theta)d\theta}$
\noindent and can be used to compute a point estimator of $\theta$. We propose to use either the maximum a posteriori (MAP) estimator $\theta_{n}$ defined as $\theta_{n} := \arg\max_{\theta\in\Theta} \mu(\theta|r_n)$ or the posterior mean estimator $\mathbf{E}(\theta|r_n): = \int_{\Theta}\theta\mu(\theta|r_n)d\theta$.
%
\subsection{Properties of the marginal posterior distribution of $\theta$}\label{sss_Properties_theta}
In this section we show two important results. The first one establishes that $\mu(\theta|r_n)$ is proportional to \eqref{eq_posterior_introduction} as stated in the introduction. To prove it we use the result of Proposition A.1 in Appendix A which shows that, once the nuisance parameter $f$ is integrated out, the posterior distribution of $\theta$ is not affected by the choice of the dominating measure $\Pi$ which only causes a transformation of the nuisance parameter.\\
\indent Therefore, we can use two different dominating measures $\Pi$: one for the definition of the nuisance parameter $f$ and one for the definition of $\mathcal{E}$, $K$ and $K^*$. More precisely, if $\sup_{x\in S}\frac{dF_*(x)}{d\Pi(x)} < \infty$ (where $\Pi$ is used to define $f$ and to construct its prior, and $F_*$ is used to define $\mathcal{E}$, $K$ and $K^*$) then, once we have specified the prior for the nuisance parameter $f=d F/d\Pi$, we deduce from it the prior of the transformation $\varphi(f):=fd\Pi/dF_*$ (see the proof of Proposition A.1 for more details). Therefore, the prior mean of $\varphi(f)$ is $f_{0\theta}d\Pi/dF_*$ and the prior covariance operator of $\varphi(f)$, denoted $\Omega_{0\theta}$, writes in terms of an o.n.b. of $L^2(S,\mathfrak{B}_{S},F_{*})$ that we still denote by $(\varphi_j)_{j\geq 0}$ and where $\{\varphi_j\}_{j=1}^d$ are equal to the moment functions $h_j(\theta,x)$, $j=1,\ldots,d$, orthonormalized with respect to $F_*$: $\Omega_{0\theta}\cdot = \sum_{j>d}\lambda_j\la\varphi_j,\cdot\ra \varphi_j$. Moreover, the operator $K$ and its adjoint $K^*$ are defined by using $ F_*$ instead of $\Pi$ so that they are operators from (resp. to) the space $\mathcal{E} = \mathcal{E}_* := L^{2}(S,\mathfrak{B}_{S},F_{*})$, and this does not change our inference on $\theta$. Hence, we can show the following theorem.
\begin{tr}\label{tr_posterior_simplified}
  Let the assumptions of Theorem A.1 in Appendix A be satisfied. Let the nuisance parameter $f$ be defined as $f:=\frac{dF}{d\Pi}$ and let $(\varphi_j)_{j\geq 0}$ be an o.n.b. of $\mathcal{E}_*:=L^2(S,\mathfrak{B}_S,F_*)$ where $\{\varphi_j(x)\}_{j=0}^d$ is an o.n.b. of $span\{1,h_1(\theta,\cdot),\ldots,h_d(\theta,\cdot)\}\subset\mathcal{E}_*$. Assume that $\sup_{x\in S}\frac{dF_*(x)}{d\Pi(x)}< \infty$. Then, $\mu(\theta|r_n)$ is proportional to the expression given in \eqref{eq_posterior_introduction} where the $\lambda_j$s do not depend on $\theta$, are positive and satisfy $\sum_{j>d}\lambda_j <\infty$.
\end{tr}
\indent This theorem, which is proved in Appendix E, shows that in order to implement our method one does not have to specify neither $r_n$ nor $K$ but only $f_{0\theta}$ (or $\varphi(f_{0\theta})$) and $\{\lambda_j\}_{j>d}$. This simplification and the definition of $K$, $K^*$ and $\mathcal{E}$ in terms of $\Pi = F_*$ must be understood in the following of the paper every time we explicitly assume $\sup_{x\in S}\frac{dF_*(x)}{d\Pi(x)} < \infty$. From a computational point of view the orthonormalization of $1$ and the moment functions $h_j(\theta,x)$, $j=1,\ldots,d$, with respect to $F_*$ can be implemented either by using the Cholesky decomposition of (the empirical counterpart of) $\mathbf{E}^*[(1,h(\theta,X)^T)^T(1,h(\theta,X)^T)]$ or by using the Gram-Schmidt process.\\
\indent The second result we are going to show\footnote{We thank Yuichi Kitamura for having suggested this research question.} establishes a link between our Bayesian procedure, the LIL procedure, the GEL estimators with quadratic criterion and the continuous-updating GMM estimator. This relationship, given in Theorem \ref{tr_equivalence_LIL} below, holds when the $\mathcal{GP}$ prior for $f|\theta$ is allowed to become diffuse. More precisely, let us rescale the prior covariance operator of $f|\theta$ by a positive scalar $c$ so that the prior of $f|\theta$ may be written, for $f_{0\theta}\in\mathcal{E}_M$, as
$$\mu(f|\theta,c) \sim \mathcal{GP}(f_{0\theta},c\Omega_{0\theta}),\quad \int h(\theta,x)f_{0\theta}(x)\Pi(dx) = 0, \quad \Omega_{0\theta}^{1/2}(1, h(\theta,\cdot)^T)^T = 0, \quad c\in\mathbb{R}_+. $$
\begin{tr}\label{tr_equivalence_LIL}
  Let the assumptions of Theorem A.1 in Appendix A be satisfied. Assume that $\sup_{x\in S}\frac{dF_*(x)}{d\Pi(x)} < \infty$, $h_j(\theta,x)\in\mathcal{R}(K^*)$ and $\varphi_l\in\mathcal{R}(K^*)$, $\forall j=1,\ldots,d$, $l>d$ and $\forall\theta\in\Theta$, and that $\mathbf{E}^*[h(\theta,X_i)h(\theta,X_i)^T]$ is nonsingular $\forall \theta\in\Theta$. Let $\mu(f|\theta,c) \sim \mathcal{GP}(f_{0\theta},c\Omega_{0\theta}),$ with $f_{0\theta}$ and $\Omega_{0\theta}$ satisfying Restrictions \ref{R1} and \ref{R2}, and $c\in\mathbb{R}_+$. Let $\mu(\theta|r_n,c)$ denote the (marginal) posterior of $\theta$ obtained by integrating out $f$ from $P^f$ with respect to $\mu(f|\theta,c)$. Then,
  \begin{displaymath}
    \lim_{c\rightarrow \infty}\mu(\theta|r_n,c) \propto \exp\left\{-\frac{1}{2}\left(\frac{1}{\sqrt{n}}\sum_{i=1}^n h(\theta,x_i)\right)^TV_n(\theta)^{-1}\left(\frac{1}{\sqrt{n}}\sum_{i=1}^n h(\theta,x_i)\right)\right\} \mu(\theta)
  \end{displaymath}
  \noindent where $V_n(\theta) := \frac{1}{n}\sum_{i=1}^nh(\theta,x_i)h(\theta,x_i)^T$.
\end{tr}
\noindent Remarks that in the theorem the limit $c\rightarrow \infty$ is taken after $f$ has been marginalized out. The result in the theorem deserves some comments. First, it shows that, as the (conditional) prior on $f$ becomes more and more diffuse, our marginal likelihood function becomes the quasi-likelihood function for GMM settings (which is the finite sample analog of the LIL) that has been considered in the literature, \textit{e.g.} by \citet{ChernozhukovHong2003} and \citet{Kim2002}. Therefore, the LIL naturally arises from a nonparametric Bayesian procedure which places a $\mathcal{GP}$ prior on the set of functions in $\mathcal{E}$ that satisfy the MRs, as the nonparametric prior becomes noninformative.\\
\indent Second, Theorem \ref{tr_equivalence_LIL} shows that, as the prior on $f$ becomes noninformative, the MAP objective function is the same (up to constants and up to the prior) as the continuous-updating GMM objective function of \cite{HansenHeatonYaron1996}. In addition, because by Theorem 2.1 in \citet{NeweySmith2004} the continuous-updating GMM estimator equals the GEL estimator with quadratic criterion, hence our MAP estimator also equals the latter.

\subsection{Properties of the MAP estimator of $\theta$}\label{s_separable_case}
By Proposition A.1 in Appendix A, our inference procedure is invariant to the choice of $\Pi$. Let us assume $\sup_{x\in S}\frac{dF_*(x)}{d\Pi(x)} < \infty$ so that $\mathcal{E} = \mathcal{E}_* := L^{2}(S,\mathfrak{B}_{S},F_{*})$. By Theorem \ref{tr_posterior_simplified} the MAP writes as:
\begin{align}
  \theta_n & = \arg\max_{\theta\in\Theta}\mu(\theta|r_n) = \arg\max_{\theta\in\Theta}\left(\log p_{n\theta}(r_n;\theta) + \log\mu(\theta)\right)\nonumber\\
  & = \arg\min_{\theta\in\Theta}\Big(\underbrace{\frac{1}{n}\sum_{i=1}^n[h(\theta,x_i)^T] \left(\E^*[h(\theta,X)h(\theta,X)^T]\right)^{-1}\frac{1}{n}\sum_{i=1}^n[h(\theta,x_i)]}_{:=A(\theta)} \nonumber \\
  &     \qquad \qquad + \sum_{j>d}\left(\frac{1}{n}\sum_{i=1}^{n}\varphi_j(x_i) - \E^*[\varphi(f_{0\theta}),\varphi_j]\right)^2\frac{1}{1 + n\lambda_j} - \frac{2}{n}\log\mu(\theta)\Big).\label{eq_posterior_separable_case}
\end{align}
Equation (\ref{eq_posterior_separable_case}) is quite useful and allows to emphasize several aspects of our methodology.
\setdefaultleftmargin{0cm}{}{}{}{}{}
    \begin{itemize}
      \item[I.] The first term in (\ref{eq_posterior_separable_case}) accounts for the MRs. Indeed, $A(\theta)$ corresponds to the orthonormalized moment function $h(\theta,x_i)$ with respect to $F_*$. Minimization of the empirical counterpart of this term corresponds to the minimization of the continuous-updating GMM criterion.
      \item[II.] The second term in \eqref{eq_posterior_separable_case} accounts for the extra information that we have, namely, the information contained in the subspace of $\mathcal{E}$ orthogonal to $span\{1,h_{1}(\theta,\cdot), \ldots, h_{d}(\theta,\cdot)\}$. This information, which is in general not exploited in MRs frameworks (in frequentist as well as in Bayesian approaches), can be exploited thanks to the prior distribution and the prior mean $f_{0\theta}$. Asymptotically, this information still plays a role if the prior is not fixed but varies with $n$ at an appropriate rate (see comment III below). On the contrary, if the prior is fixed then, as $n\rightarrow \infty$, the second term of \eqref{eq_posterior_separable_case} converges to $0$ since $(1 + n\lambda_j)^{-1}\rightarrow 0$. 
      \item[III.] Expression (\ref{eq_posterior_separable_case}) makes an explicit connection between the parametric case (infinite number of MRs) and the semiparametric case (where only the first $d$ MRs hold). The semiparametric case corresponds to the classical GEL or GMM approach while the parametric case corresponds to the maximum likelihood estimator (MLE). Indeed, the prior distribution for $f$ specifies a parametric model for $f_{0\theta}$ which satisfies the $d$ MRs and eventually other ``extra'' MRs. The eigenvalues $\lambda_{j}$ of the prior covariance operator play the role of weights of the ``extra'' MRs and represent our ``beliefs'' concerning these restrictions. When we are very confident about these ``extra'' conditions, or equivalently we believe that $f_{0\theta}$ is close to $f_{*}$, then the $\lambda_{j}$s are close to zero or converge to $0$ faster than $n^{-1}$ as $n\rightarrow \infty$. So, the prior distribution for $f$ is degenerate on $f_{0\theta}$ (as $n$ increases). In that case, the MAP estimator will essentially be equivalent to the MLE that we would obtain if we use the prior mean function $f_{0\theta}$ as the likelihood. When we are very uncertain about $f_{0\theta}$ then the $\lambda_{j}$s are very large and may tend to $+\infty$ (uninformative prior). In this case the MAP estimator will be close to the continuous-updating GMM estimator (up to a prior on $\theta$).
    \end{itemize}

\section{Asymptotic Analysis}\label{s_asymptotic_analysis}
In this section we focus on the frequentist asymptotic properties of our approach for $n\rightarrow\infty$. For this analysis we use the true probability measure $P^*$ which corresponds to the true DGP $F_*$ and maintain the assumption that $\sup_{x\in S}\frac{dF_*(x)}{d\Pi(x)}<\infty$ where $\Pi$ is used to define the nuisance parameter $f$ and $F_*$ to define $K$, $K^*$ and $\mathcal{E}$ (see the discussion below Proposition A.1). We show: \textit{(i)} consistency of the posterior of $\theta$ (Theorem \ref{tr_posterior_consistency_posterior_theta}), \textit{(ii)} frequentist consistency of the MAP estimator $\theta_n$ (Theorem \ref{tr_posterior_consistency_theta}), and \textit{(iii)} convergence in Total Variation of $\mu(\theta|r_n)$ towards a normal distribution (section \ref{ss_asymptotic_normality}).

\setdefaultleftmargin{0.5cm}{}{}{}{}{}
\subsection{Posterior Consistency}\label{ss_posterior_consistency}
We first state the following assumption.
\begin{assump}
  \textit{(1) the true parameter $\theta_*$ belongs to the interior of a compact convex subset $\Theta$ of $\mathbb{R}^p$ and is the unique solution of $\mathbf{E}^{*}[h(\theta,X)] = 0$;\label{Ass_A0}
  (2) the singular functions $\{\varphi_j(X)\}_{j>d}$ and $h(\theta,X)$ are twice continuously differentiable in a neighborhood of $\theta_*$ with probability $1$; (3) the prior distribution $\mu(\theta)$ is a continuous probability measure that admits a density with respect to the Lebesgue measure and is positive on a neighborhood of $\theta_*$; (4) the function $h(\theta,X)$ is continuous at each $\theta\in\Theta$ with probability $1$, and there is $\delta(X)$ with $\|h(\theta,X)\|\leq \delta(X)$ for all $\theta\in\Theta$ and $\mathbf{E}[d(X)]< \infty$.}
\end{assump}
\noindent Assumption \ref{Ass_A0} is a standard assumption in the literature on MRs. The next theorem gives concentration of the posterior distribution around the true value $\theta_*$.
\begin{tr}\label{tr_posterior_consistency_posterior_theta}
  Let the assumptions of Theorem \ref{tr_posterior_simplified} and Assumption \ref{Ass_A0} be satisfied and assume $\sup_{x\in S}\frac{dF_*(x)}{d\Pi(x)} < \infty$. Then, for any sequence $M_n\rightarrow \infty$,
    \begin{equation}\label{eq_concentration_posterior}
      \mu\left( \sqrt{n}\|\theta - \theta_*\|>M_n|r_n\right) \rightarrow 0\qquad \textrm{in $P^*$-probability as $n\rightarrow \infty$.}
    \end{equation}
\end{tr}
Given the result of the theorem, the next result follows.
\begin{tr}\label{tr_posterior_consistency_theta}
  Let the assumptions of Theorem \ref{tr_posterior_consistency_posterior_theta} be satisfied. Then, $\theta_n\xrightarrow{p} \theta_*$ in $P^*$-probability as $n\rightarrow \infty$.
\end{tr}

\subsection{Asymptotic Normality}\label{ss_asymptotic_normality}
In this section we establish asymptotic normality of $\mu(\theta|r_n)$. This result applies to the case $d>p$ (which is our main interest) as well as to the case $d = p$. In appendix H we establish, under different assumptions, asymptotic normality of $\mu(\theta|r_n)$ for the case where $d = p$ and an alternative prior for $\theta$ deduced from the prior for $f$.\\
\indent For some $\tau\in \mathbb{R}^p$, let $s_n(\tau) := p_{n,\theta_*+\tau/\sqrt{n}}(r_n;\theta_* + \tau/\sqrt{n})$, where $p_{n,\theta}$ is defined in Theorem A.1. We assume that there exist a random vector $\tilde{\ell}_*$ and a nonsingular matrix  $\tilde{I}_{*}$ (that depend on the true $\theta_*$ and $f_*$) such that the sequence $\tilde{\ell}_*$ is bounded in probability, and satisfy
\begin{equation}\label{eq_Integral_LAN}
  \log\frac{s_n(\tau)}{s_n(0)} = \frac{1}{\sqrt{n}}\tau^T\widetilde{I}_*\widetilde{\ell}_* - \frac{1}{2}\tau^T\widetilde{I}_{*}\tau + o_p(1)
\end{equation}
\noindent for every random sequence $\tau$ which is bounded in $P^*$-probability. Condition \eqref{eq_Integral_LAN} is known as the integral local asymptotic normality assumption which is used to prove asymptotic normality of semiparametric Bayes procedures, see e.g. \citet{BickelKleijn2012}. In Appendix G.2 we prove that, if $\sup_{x\in S}\frac{dF_*(x)}{d\Pi(x)} < \infty$ and other assumptions given in the appendix hold, then equation \eqref{eq_Integral_LAN} holds with
\begin{equation}\label{eq_post_var}
  \widetilde{I}_* = \mathbf{E}^*\left[\frac{\partial h(\theta_*,x)^T}{\partial \theta}\right]\left(\mathbf{E}^*[h(\theta_*,X)h(\theta_*,X)^T]\right)^{-1}\mathbf{E}^*\left[\frac{\partial h(\theta_*,X)}{\partial \theta^T}\right]
\end{equation}
\noindent and $\frac{1}{\sqrt{n}}\widetilde I_*\widetilde \ell_* = -\mathbf{E}^*\left[\frac{\partial h(\theta_*,x)}{\partial  \theta^T}\right]\left(\mathbf{E}^*[h(\theta_*,x)h(\theta_*,x)^T]\right)^{-1}\frac{1}{\sqrt{n}}\sum_{i=1}^n h(\theta_*,x)$ if $\mathbf{E}^*\left[h(\theta_*,X)h(\theta_*,X)^T\right]$ is nonsingular. For two probability measures $P_1$ and $P_2$ absolutely continuous with respect to a positive measure $Q$, define the total variation (TV) distance as $||P_1-P_2||_{TV} := \int|f_{1} - f_{2}|dQ$ where $f_{1}$ and $f_{2}$ are the Radon-Nikodym derivatives of $P_1$ and $P_2$, respectively, with respect to $Q$. The following theorem shows that under \eqref{eq_Integral_LAN} the posterior distribution of $\sqrt{n}(\theta - \theta_*)$ converges in the TV distance to a Normal distribution with mean $\Delta_* := \frac{1}{\sqrt{n}}\widetilde{\ell}_*$ and variance $\widetilde{I}_*^{-1}$. From \eqref{eq_post_var} we see that the posterior distribution has the same asymptotic variance as the efficient GMM estimator.
\begin{tr}\label{tr_TV_overidentified}
  Let the assumptions of Theorem \ref{tr_posterior_simplified} and Assumption \ref{Ass_A0} be satisfied. Assume that \eqref{eq_concentration_posterior} and \eqref{eq_Integral_LAN} hold. If $\mu(\sqrt{n}(\theta - \theta_*)|r_n)$ denotes the posterior of $\sqrt{n}(\theta - \theta_*)$, then:
  \begin{equation}\label{eq_TV_overidentified}
    \|\mu(\sqrt{n}(\theta - \theta_*)|r_n) - \mathcal{N}(\Delta_*,\widetilde{I}_{*}^{-1})\|_{TV}\rightarrow 0\qquad \textrm{in $P^*$-probability as $n\rightarrow \infty$.}
  \end{equation}
\end{tr}
%
\section{Implementation}\label{s_implementation}
In this section we first explain the numerical implementation of our procedure and then show the results of two simulations and an application to demand estimation for airline seats. To implement our procedure one has to use the expression of the marginal posterior distribution of $\theta$ as given in \eqref{eq_posterior_introduction}. To implement \eqref{eq_posterior_introduction} one has to replace $F_*$ by the empirical cdf $F_n$ and truncate the infinite sum at some $J<n$. Then, a Metropolis-Hastings (M-H) algorithm has to be used to draw from $\mu(\theta|r_n)$. Before starting the M-H algorithm, one has to: (1) choose $J$, (2) choose the sequence $\lambda_j$ of eigenvalues of $\Omega_{0\theta}$ such that $\sum_j \lambda_j < \infty$ (for instance: $\lambda_j = j^{-\alpha}$ with $\alpha > 1$, or $\alpha^{-j}$ with $\alpha>1$). For simplicity the prior on $\theta$ is chosen to be flat. In the simulations and empirical application we have used a random walk M-H algorithm to draw from $\mu(\theta|r_n)$ which consists of the following steps (other M-H algorithms can be used).
\paragraph{Algorithm 1: (random walk Metropolis-Hastings).}
    \begin{enumerate}
      \item conditionally on $\theta^{(j-1)}$, generate $\theta^{(j)}$ from the following auxiliary distribution: $\theta^{(j)}\sim\mathcal{N}(\theta^{(j-1)},\xi D)$ where $D$ is equal to the estimated asymptotic posterior variance given in Theorem \eqref{tr_TV_overidentified} and $\xi>0$ is tuned to guarantee an appropriate acceptance rate 
          as suggested in \cite{GRG1996};
      \item evaluate the $d$-vector function $h(\theta,x)$ at $\theta^{(j)}$ and at the $n$ observations $x_{1},\ldots,x_{n}$, and orthonormalize them by computing
          $$A_n(\theta^{(j)}) := \frac{1}{n}\sum_{i=1}^n[h(\theta^{(j)},x_i)^T] \left(\frac{1}{n}\sum_{i=1}^n[h(\theta^{(j)},x_i)h(\theta^{(j)},x_i)^T]\right)^{-1}\frac{1}{n}\sum_{i=1}^n[h(\theta^{(j)},x_i)];$$
      \item \emph{construction of $(\varphi_j)_{d+2\leq j\leq J}$:} (a) randomly generate $J-(d+1)$ independent vectors $v_j$ of dimension $n$ for $j= d+2, \ldots, J$; (b) orthonormalize the set of vectors $\{\iota,(h_1(\theta^{(j)},x_1),\ldots,h_1(\theta^{(j)},x_n))^T, \ldots, (h_d(\theta^{(j)},x_1),\ldots,h_d(\theta^{(j)},x_n))^T, (v_j)_{d+2\leq j\leq J}\}$, where $\iota$ is a $n$-vector of ones, by using \textit{e.g.} the Gram-Schmidt process; (c) set $(\varphi_j)_{d+2\leq j\leq J}$ equal to the last $J-(d+1)$ of these orthonormalized vectors and denote them by $(\varphi_k(\theta^{(j)},x_1),\ldots,\varphi_k(\theta^{(j)},x_n))_{d+2\leq k\leq J}$;
      \item for $\theta^{(j-1)}$ compute the same quantities as in 2 and 3;
      \item for both $\theta = \theta^{(j-1)}$ and $\theta = \theta^{(j)}$ construct a prior mean function $f_{0\theta}$. This can be done either by using the method described in section \ref{ss_prior_distribution} if the MRs are multivariate polynomial functions, or by taking a $f_{0\theta}$ which is orthogonal to the moment functions $h(\theta,x)$ and to the $(J - (d+1))$ elements of the basis $(\varphi_j)_{d+2\leq j\leq J}$ with respect to $F_n$;
      \item compute the logarithm of the acceptance probability of the M-H sampler as follows:
      \footnotesize{
      \begin{flalign*}
        &\alpha(\theta^{(j-1)},\theta^{(j)}) = - \frac{n}{2}\left[A_n(\theta^{(j)}) + \sum_{k>d+1}\frac{\left(\frac{1}{n}\sum_{i=1}^{n}\varphi_k(\theta^{(j)},x_i) - \frac{1}{n}\sum_{i=1}^n[f_{0\theta^{(j)}}(x_i)\varphi_k(\theta^{(j)},x_i)]\right)^2}{1 + n\lambda_k}\right]&\\
        &+ \frac{n}{2}\left[A_n(\theta^{(j-1)}) - \sum_{k>d+1}\frac{\left(\frac{1}{n}\sum_{i=1}^{n}\varphi_k(\theta^{(j-1)},x_i) - \frac{1}{n}\sum_{i=1}^n[f_{0\theta^{(j-1)}}(x_i)\varphi_k(\theta^{(j-1)},x_i)]\right)^2}{1 + n\lambda_k}\right];&
      \end{flalign*}
      }
      \normalsize
      \item draw $U\sim Uniform[0,1]$ and accept $\theta^{(j)}$ if $\log(U) \leq \log \alpha(\theta^{(j-1)},\theta^{(j)})$, otherwise set $\theta^{(j)} = \theta^{(j-1)}$;
      \item iterate steps 1-7 and discard the first draws corresponding to the burn-in period.
    \end{enumerate}
%
\subsection{Overidentified case}\label{ss_overidentified_case}
Let us consider the case in which $x$ is univariate and the one-dimensional parameter of interest $\theta$ is characterized by the moment conditions $\mathbf{E}^{F}[h(\theta,X)] = 0$ with $h(\theta,x) = (x - \theta, 2\theta^{2} - x^2)^T$. These moment conditions are satisfied for instance for $F$ an exponential distribution with parameter $\theta$. The prior $\mu(\theta)$ is uniform on $[\theta_* - 1,\theta_*+1]$.\\
\indent In our simulation, we generate $n$ observations $x_1,\ldots,x_n$ from an exponential distribution with parameter $\theta_* = 2$. We use polynomially decreasing eigenvalues for $\Omega_{0\theta}$: $\lambda_j = j^{- 1.7}$ and use only the first $J = 57$ eigenvalues. Remark that truncating the series is equivalent to assume that $\lambda_j=\infty$ for $j>J$, that is, the prior is diffuse only on some directions (the ones corresponding to $j>J$). For $n$ large, truncation does not affect the results much. On the contrary, for $n$ small, truncation of the series does affect the results. This can be seen by comparing the results for $n = 100$ and $n=500$ in Tables 1 and 6. The prior mean function $f_{0\theta}$ is set equal to a $\mathcal{N}(\theta,\theta^2)$.\\
\indent To draw from the posterior distribution of $\theta$, we use the random walk M-H algorithm described in \textbf{Algorithm 1} with auxiliary distribution a $\mathcal{N}(\theta^{(j-1)},\xi D)$ distribution, where $D$ is equal to the estimated asymptotic posterior variance, that is, $$D = \left.\frac{1}{n}\sum_{i=1}^n\frac{\partial h(\theta,x_1)'}{\partial \theta} \left[\frac{1}{n}\sum_{i=1}^nh(\theta,x_i) h(\theta,x_i)'\right]^{-1} \frac{1}{n}\sum_{i=1}^n\frac{\partial h(\theta,x_i)'}{\partial \theta}\right|_{\theta = \theta_n},$$ and $\theta_n$ is the MAP estimator. We have tried different scenarios involving different sample size $n$, different $\xi$ and different initial values $\theta^{(0)}$ for the chain. In Table \ref{Table_Simulation_1_Overidentified} and Figures 1-6 we report the results for all these scenarios. In Table \ref{Table_Simulation_1_Overidentified} we report the posterior mean, standard deviation, median, MAP, the lower and upper bounds of the $95\%$ HPD-credible region, and the Acceptance Rate. We notice that the posterior standard deviation for the case with $n=2000$ is about half the posterior standard deviation with $n=500$ as predicted by the Bernstein-von Mises theorem. Figures 1-6 in Appendix C.1 show, together with the histogram of the draws from the posterior of $\theta$, three plots that are used for diagnose the convergence of our MCMC: the trace plot, the autocorrelation function and the running mean plot. All these plots indicates that the series has converged. The histogram, where we have superposed the Gaussian distribution shows how the posterior distribution can be approximated by a Gaussian distribution. The frequentist estimator of $\theta$, denoted by $\bar{X}$ and computed as the sample mean, is also reported in Table \ref{Table_Simulation_1_Overidentified} for different sample sizes.\\
\indent For comparison, we show in Table 5 in Appendix C.1 the results obtained with the method proposed by \cite{ChernozhukovHong2003}. The corresponding graphs are in Figures 7-9 in Appendix C.1. The difference between the two methods is evident for $n=100$ but almost disappears for $n=500$. Indeed, increasing the sample size plays the same role as making the prior distribution more and more diffuse. Finally, Table 6 in Appendix C.1 shows the impact of the truncation on the inference for $\theta$. Here, $J$ is fixed equal to $100$ instead of $57$ as before which means that the effect of the prior is more important since we keep more terms in the second sum in \eqref{eq_posterior_introduction}. Indeed, we see that for small sample size ($n=100$) the bias is more important than in the case with $J=57$. However as long as the sample size increases from $n=100$ to $n=500$ the bias decreases.

\begin{table}[h]
\centering
\footnotesize{
  \begin{tabular}{c||ccccccc}
    \hline
    \hline
    & Mean & SD & Median & MAP & Lower & Upper & AR \\
    \hline
    \multirow{4}{*}{$n = 100$} & &\multicolumn{5}{c}{\footnotesize{$\theta^{(0)} = \bar{X} + \mathcal{N}(0,1) = 1.3726$, $\xi = 0.31$}} & \\
    & 1.8165 & 0.1385 & 1.8200 & 1.8463  & 1.5361 & 2.0746 & $44.89\%$\\
    \cline{2-8}
    & &\multicolumn{4}{c}{\footnotesize{$\theta^{(0)} = \bar{X} = 1.8963$, $\xi = 0.28$}}&\\
     & 1.8168 & 0.1407 & 1.8189 & 1.8463 & 1.5296 & 2.0806 & $45.05\%$\\
    \hline
    \multirow{4}{*}{$n = 500$} & &\multicolumn{5}{c}{\footnotesize{$\theta^{(0)} = \bar{X} + \mathcal{N}(0,1) = 1.2694$, $\xi = 0.16$}} & \\
    & 1.9995 & 0.0891 & 1.9991 & 2.0032 & 1.8261 & 2.1758 & $42.19\%$\\
    \cline{2-8}
    & &\multicolumn{5}{c}{\footnotesize{$\theta^{(0)} = \bar{X} = 2.0032$, $\xi = 0.18$}}&\\
     & 1.9995 & 0.0905 & 1.9987 & 2.0032 & 1.8247 & 2.1794 & $44.67\%$\\
    \hline
    \multirow{4}{*}{$n = 2000$} & &\multicolumn{5}{c}{\footnotesize{$\theta^{(0)} = \bar{X} + \mathcal{N}(0,1) = 3.2254$, $\xi = 0.12$}} & \\
    & 2.0061 & 0.0439 & 2.0057 & 2.0075 & 1.9207 & 2.0939 & $29.46\%$\\
    \cline{2-8}
    & &\multicolumn{5}{c}{\footnotesize{$\theta^{(0)} = \bar{X} = 2.0075$, $\xi = 0.12$}}&\\
    & 2.0067 & 0.0436 & 2.0068 & 2.0075 & 1.9202 & 2.0926 & 29.03\%\\
    \hline
    \hline
  \end{tabular}}
\caption{{\footnotesize Results obtained with our Gaussian Process approach for $J=57$, different sample sizes and different settings by using an MCMC with $20,000$ draws after a burn-in of $10,000$. True value of theta is $2$. Mean = posterior mean, SD = posterior standard deviation, Median = posterior median, MAP = maximum a posteriori, Lower = lower bound (resp. Upper = upper bound) of the $95\%$ HPD-credible region, AR = acceptance rate.}}\label{Table_Simulation_1_Overidentified}
\end{table}

\subsection{Linear instrumental regression estimation}\label{ss_IV_simulation}
In this section we study the finite sample properties of our procedure in the setting of a linear regression model where some covariates are endogenous. The parameters of interest are identified through instrumental variables. The data generating process is the following: $Y = \theta_0 W_0 + \theta_1^T W_1 + 2U$, $W := (W_0, W_1^T)^T$, $W_0 = \alpha_1 \tilde{Z} + \alpha_2^T W_1 + V$, $Z := (\tilde{Z},W_1^T)^T$, $\tilde{Z}$ is $1$-dimensional and $W_1$ has dimension $2$. Moreover, $(U,V)\sim \mathcal{N}(0,\left(\begin{array}
    {cc}1 & 0.6\\
    0.6 & 1  \end{array}\right))$, and $Z\sim\mathcal{N}(0,\Sigma)$ independent of $(U,V)$ with
\noindent $\Sigma = (0.5^{|i-j|})_{ij}$. The true value of the parameters is: $\theta = (\theta_0,\theta_1^T)^T = (2,1.5,-3)^T$ and $\alpha = (\alpha_1,\alpha_2^T)^T = (10,-2,1)^T$. The moment functions $h(\theta,X)$ are given by $h(\theta, X)=[(Y - \theta^T W)W_1^T,(Y - \theta^T W)\tilde{Z}]^T$
\noindent for $X = (Y,W_0,Z^T)^T$. We orthonormalize the functions $\{h_1(\theta,X),h_2(\theta,X),h_3(\theta,X)\}$ by using the empirical measure and then complete this basis to have a total of $57$ basis functions. The eigenvalues of the prior covariance operator are set equal to $\lambda_j = j^{-1.7}$ and the prior mean function $f_{0\theta}$ is constructed by taking a $\mathcal{N}_3(0,\Sigma)$ with $\Sigma$ restricted to satisfy the MRs for a given $\theta$. The prior for $\theta$ is specified as a product of uniform distributions. To draw from the posterior distribution of $\theta$, we use a random walk M-H algorithm given in \textbf{Algorithm 1} with auxiliary distribution a $3$-dimensional $\mathcal{N}(\theta^{(i-1)},\xi D)$ distribution, where $\theta^{(i-1)}$ denotes the value of $\theta$ at the previous iteration, $\xi>0$, and $D$ is equal to the estimated asymptotic posterior variance-covariance matrix. We have tried different scenarios involving different sample size $n$ and different $\xi$. The results are given in Table \ref{Table:estimators:IV:Simulation} where we report the posterior mean, standard deviation and  median, the MAP, the $95\%$ HPD-credible region and the two stage least square estimate $\widehat\theta_{2SLS}$. Figures 10-12 in Appendix C.2 show for the three sample sizes considered the histograms of the draws from the posteriors of $\theta$ and plots that are used for diagnose the convergence of our MCMC: the trace plots, the autocorrelation function and the running mean plots. For comparison, we report in Table 7 in Appendix C.2 the results obtained with the method of \cite{ChernozhukovHong2003}. We see that the results obtained with the two methods are quite close. This is due to the fact we have used a small value for the truncation parameter $J$, see the discussion in Section \ref{ss_overidentified_case}.
\begin{table}[h]
\centering
\footnotesize{
  \begin{tabular}{c||c||c||cccccccc}
    \hline
    \hline
    &  & & Mean & SD & Median & MAP & Lower & Upper & $\widehat{\theta}_{2SLS}$\\
    \hline
    \multirow{3}{*}{$n = 100$} & $\xi = 0.39$& $\theta_0$ & 2.0066 & 0.0283 & 2.0058 & 2.0054 & 1.9531 & 2.0658 & 2.0054 \\
    &  \multirow{3}{*}{$AR = 37.19\%$} & $\theta_{11}$ & 1.1130 & 0.3349 & 1.1220 & 1.1746 & 0.4174 & 1.7609 & 1.1745\\
    & & $\theta_{12}$ & -2.9707 & 0.2846 & -2.9817 & -2.9980& -3.4880 & -2.3647 & - 2.9980\\
    \hline
    \multirow{3}{*}{$n = 500$} & $\xi = 0.244$ & $\theta_0$ & 1.9980 & 0.0109 & 1.9979 & 1.9986 & 1.9759 & 2.0193 & 1.9988\\
    & \multirow{3}{*}{$AR = 45.83 \%$} & $\theta_{11}$ & 1.4666 & 0.1215 & 1.4665 & 1.4786 & 1.2308 & 1.7015 & 1.4712\\
    & & $\theta_{12}$ & -2.9629 & 0.1074 & -2.9642 & -2.9528 & -3.1687 & -2.7439 & -2.9680\\
    \hline
    \multirow{3}{*}{$n = 1000$} & $\xi = 0.26$ & $\theta_0$ & 2.0025 & 0.0076 & 2.0024 & 2.0034 & 1.9876 & 2.0179 & 2.0027\\
    & \multirow{3}{*}{$AR = 43.58\%$} & $\theta_{11}$ & 1.4190 & 0.0845 & 1.4196 & 1.4177 & 1.2481 & 1.5835 & 1.4229\\
    & & $\theta_{12}$ & -2.9875 & 0.0760 & -2.9867 & -2.9807 & -3.1396 & -2.8360 & -2.9901\\
    \hline
    \hline
  \end{tabular}}
\caption{{\footnotesize Linear instrumental regression. True value of $\theta$ is $\theta_0 = 2$, $\theta_{11} = 1.5$ and $\theta_{12} = -3$. Results obtained with our Gaussian Process approach for $J = 57$, $\theta^{(0)} = \widehat{\theta}_{2SLS}$ and different sample sizes by using an MCMC with $20,000$ draws after a burn-in of $10,000$. Mean = posterior mean, SD = posterior standard deviation, Median = posterior median, MAP = maximum a posteriori, Lower = lower bound (resp. Upper = upper bound) of the $95\%$ HPD-credible region, AR = acceptance rate.}}\label{Table:estimators:IV:Simulation}
\end{table}

\subsection{Application to demand estimation}\label{ss_Application_demand_estimation}
In this application we consider estimation of a demand function for airline seats that we model as $\log (passen) = \theta_0 + \theta_1 \log(fare) + \theta_2 \log(dist) + \theta_3 (\log(dist))^2 + \epsilon,$ \noindent where $passen$ denotes the average passengers per day, $fare$ is the average airfare and $dist$ is the route distance. This is an example of simultaneous equations model where demanded quantity ($passen$) and price ($fare$) are simultaneously determined. Therefore, we need an additional instrumental variable to consistently estimate the demand equation. Following \cite[Chapter 16]{Wooldridge2013}, we use a measure of market concentration as an instrument, namely, we use the share of business accounted for by the largest carrier, denoted by $concen$ in the following. This is a valid instrument under the assumption that passenger demand is determined only by air fare, so that, once price is controlled for, passengers are indifferent about the share of business accounted for by the largest carrier. Therefore, the parameter $\theta := (\theta_0,\theta_1,\theta_2,\theta_3)'$ is identified by the following set of moment conditions:
\begin{equation}
  \mathbf{E}[\epsilon(1,concen,\log(dist),(\log(dist))^2)^T] = 0.
\end{equation}
\noindent To estimate this model we use the data set ``AIRFARE'' provided by \cite{Wooldridge2013}\footnote{Available at \url{http://www.cengage.com/aise/economics/wooldridge_3e_datasets/}
. } and extract from it the year $1997$. This data set concerns airline companies on routes in the United States. The sample size is of $1149$ observations.\\
\indent We orthonormalize the functions $\{\epsilon, \epsilon *concen, \epsilon*\log(dist),\epsilon*(\log(dist))^2\}$ as elements of $L^2(\mathbb{R}^5,\mathfrak{B},F_*)$ with $F_*$ replaced by the empirical measure $F_n$ and complete the basis to have a total of $57$ basis functions. The prior covariance $\Omega_{0\theta}$ is constructed by using these basis functions and polynomially decreasing eigenvalues $\lambda_j = j^{-1.7}$. The prior mean function $f_{0\theta}$ is constructed by taking a multivariate normal distribution $\mathcal{N}_5(m,\Sigma)$, where $m = (\theta_0,0,0,0,0)^T$ and $\Sigma$ is restricted to satisfy the MRs for a given $\theta$. We use a product of uniform distributions as a prior for the elements of $\theta$.\\
\indent We draw from the posterior distribution of $\theta$ by using a random walk M-H algorithm with the auxiliary distribution given in \textbf{Algorithm 1}. Table \ref{Table:estimators:Airfare:Exact_Identification} reports the posterior mean, standard deviation and median, the MAP, the $95\%$ HPD-credible region and the acceptance rate obtained with our $\mathcal{GP}$-approach. For comparison, we also provide in this table the two stage least squares estimate $\widehat{\theta}_{2SLS}$ and the results obtained with the method of \cite{ChernozhukovHong2003}. The posterior mean estimate of the price elasticity $\theta_1$ is around $-1.38$ which shows that, as expected, passenger demand drops of about $1.38\%$ due to a one percent increase in fare. Figure 13 in Appendix C.3 shows the histograms of the draws from the posteriors of $\theta$ and plots that are used for diagnose the convergence of our MCMC: the trace plots, the autocorrelation function and the running mean plots.\\
\begin{table}[h]
\centering
\footnotesize{
  \begin{tabular}{c||c||c||cccccccc}
    \hline
    \hline
    &  & & Mean & SD & Median & MAP & Lower & Upper & $\widehat{\theta}_{2SLS}$\\
    \hline
    \multirow{4}{*}{$\mathcal{GP}$} & $\xi = 1.2$, $J=57$& $\theta_0$ & 19.4621 & 4.2176 & 19.1596 & 18.0691 & 12.1774 & 28.8636 & 18.0138 \\
    &  \multirow{3}{*}{$AR = 44.19\%$} & $\theta_{1}$ & -1.3763 & 0.5573 & -1.3162 & -1.1788 & -2.6836 & -0.4653 & -1.1740\\
    & & $\theta_{2}$ & -2.3960 & 0.8705 & -2.3643 & -2.1872 & -4.2297 & -0.7981 & -2.1757\\
    & & $\theta_3$ & 0.2102 & 0.0760 & 0.2059 & 0.1880 & 0.0740 & 0.3759 & 0.1870\\
    \hline
    \multirow{4}{*}{CH} & $\xi = 1.2$& $\theta_0$ & 19.5730 & 4.4161 & 19.1739 & 18.0979 & 12.2220 & 29.0699 & 18.0138 \\
    &  \multirow{3}{*}{$AR = 43.71\%$} & $\theta_{1}$ & -1.3865 & 0.5692 & -1.3257 & -1.1815 & -2.6670 & -0.4732 & -1.1740\\
    & & $\theta_{2}$ & -2.4182 & 0.9045 & -2.3583 & -2.1930 & -4.3680 & -0.8340 & -2.1757\\
    & & $\theta_3$ & 0.2122 & 0.0796 & 0.2058 & 0.1886 & 0.0775 & 0.3838 & 0.1870\\
    \hline
    \hline
  \end{tabular}}
\caption{{\footnotesize Demand estimation. $\mathcal{GP}$ (resp. CH) means that the results are obtained with our Gaussian Process (resp. \cite{ChernozhukovHong2003}) approach. We use an MCMC with $50,000$ draws after a burn-in of $50,000$ and initial value $\theta^{(0)} = \widehat{\theta}_{2SLS}$. Mean = posterior mean, SD = posterior standard deviation, Median = posterior median, MAP = maximum a posteriori, Lower = lower bound (resp. Upper = upper bound) of the $95\%$ HPD-credible region, AR = acceptance rate.}}\label{Table:estimators:Airfare:Exact_Identification}
\end{table}
\indent Next, one may wonder whether adding a power of the $concen$ variable can provide a valid instrument. To test this, we implement our testing procedure described in Appendix B. We want to test the validity of the overidentifying MR $\mathbf{E}[\epsilon *concen^2]=0$, maintaining the previous MRs, so we are in situation (B) of Appendix B. We specify a prior on $\lambda_5$ as $\mu(\lambda_{5}) = 0.5 \1\{\lambda_{5} = 0\} + 0.5 g_{\lambda_5}(\lambda_{5})\1\{\lambda_5> 0\}$, where $g_{\lambda_5} = Beta(2,2)$. So, we are comparing two models: Model 1, where the overidentifying MR does not hold, \textit{i.e.} $\mathbf{E}[\epsilon *concen^2]\neq0$ and $\lambda_5 > 0$, and Model 2, where the overidentifying MR holds, \textit{i.e.} $\mathbf{E}[\epsilon * concen^2] = 0$ and $\lambda_5 = 0$. The null hypothesis is $\mathcal{H}_0: \lambda_5 = 0$ and we test it by using the Bayes factor in favour of Model 1:
$B_{12} = \frac{0.5\int_{0}^1\int_{\Theta}p_{n\theta}(r_n;\theta,\lambda_5)\mu(\theta)d\theta g_1(\lambda_5)d\lambda_5}{0.5\int_{\Theta}p_{n\theta}(r_n;\theta,\lambda_5 = 0)\mu(\theta)d\theta}.$
Therefore, Model 1 is selected if $B_{12}>1$. By using the data set described above, we obtain $\log(B_{12}) = 9.5454$ which is a strong evidence against the validity of the overidentifying MR. We have also run a regression of $lfare$ on $(concen,concen^2,ldist,ldist^2)$ and we find that $concen^2$ is statistically significant at the level of $10\%$ but not at $5\%$. This, together with the result of our test, suggests that $concen^2$ is not a valid instrument.

\section{Conclusions}\label{s_conclusions}
This paper develops Bayesian inference for econometric models characterized via MRs. While Bayesian inference has many advantages, it is challenging when the model is characterized only by MRs so that the data distribution is left unrestricted apart from these MRs. The first difficulty is due to the fact that a likelihood function is not available. A second difficulty arises when one wants to impose the MRs in the nonparametric prior for the data distribution.\\
\indent The main contribution of this paper is to construct a prior on the nonparametric data distribution that allows to incorporate the MRs. This prior is based on a Gaussian process. With this prior and the convenient sampling model that we use, we are able to obtain a marginal posterior distribution $\mu(\theta|r_n)$ for the parameter of interest $\theta$ that is particularly convenient. We show that implementation of our procedure only requires few steps.\\
\indent Remarkably, we show that $\mu(\theta|r_n)$ is proportional to a monotonic transformation of the continuous-updating GMM criterion. We also show that the LIL-based procedure, which has been often used in the past econometric literature due to its computational convenience but which lacks a rigorous derivation, arises from our degenerate $\mathcal{GP}$ prior as it becomes noninformative.

\paragraph{Supplementary Material:} it contains: (I) a testing procedure, (II) further details about the simulations and empirical application, (III) all the proofs of the results in the paper, and (IV) an alternative procedure to be used when the model is exactly identified.

\appendix

\bibliographystyle{alphanat}
\bibliography{AnnaBib}

\end{document}